\documentclass[12pt]{article}
\usepackage[cp1251]{inputenc}
\usepackage[T2A]{fontenc}
\usepackage[english]{babel}
\usepackage{amsthm}
\usepackage{epsfig, graphicx, xcolor}
\usepackage{amssymb}
\usepackage{amsmath, amsfonts}
\newtheorem{theorem}{Theorem}
\newcommand{\RR}{\mathbb{R}}

\newcommand{\ZZ}{\mathbb{Z}}

\newtheorem{lemma}{Lemma}

\begin{document}

\title{On One Problem of Optimization of Approximate Integration
\thanks{This paper was published in Russian in Studies in Modern Problems of Summation and Approximation
of Functions and Their Applications, Collect. Sci. Works, Dnepropetrovsk, 1984, pp. 3 -13}}

\author{V. F. Babenko}
\date{}

\maketitle

\begin{abstract}
It is proved that interval quadrature formula of the form
$$
q(f)=\sum\limits_{k=1}^nc_k\frac
1{2h}\int\limits_{x_k-h}^{x_k+h}f(t)dt
$$
($c_k\in \RR , \, x_1+h<x_2-h<x_2+h<...<x_n-h<x_n+h<x_1+2\pi -h$)
with equal $c_k$ and equidistant $x_k$ is optimal among all such
formulas for the class $K*F_1$ of convolutions of a $CVD$-kernel $K$
with functions from the unite ball of the space $L_1$ of
$2\pi$-periodic integrable functions.

\medskip

\end{abstract}

1. Let $C$ and $L_p\; (1\le p\le\infty )$ be the spaces of $2\pi$
-periodic functions endowed with corresponding norms; $\|\cdot \|_p$
-- norm in $L_p$. The convolution of functions $K\in L_1$ (kernel of
convolution) and $\phi\in L_1$ is defined by equality
$$
K*\phi (x)=\int\limits_0^{2\pi}K(x-t)\phi (t)dt.
$$

Given a kernel $K$ set $\mu =\mu (K)=1$, if
$\int\limits_0^{2\pi}K(t)dt=0,$ and $\mu =\mu (K)=0$, if
$\int\limits_0^{2\pi}K(t)dt\neq 0.$ Denote by $\nu (f)$ the number
of sign changes over a period of a $2\pi$-periodic function $f$.

A kernel $K$ is called a $CVD$-kernel (denoted by $K\in CVD$), if for any function of the form
$$
f(x)=a\mu +K*\phi (x)\eqno (1)
$$
($a\in\RR,\, \phi\in C,\, \phi\perp \mu$) the inequality $\nu (f)\le\nu (\phi)$ holds.

Let $\psi (x)$ be an entire function of the form
$$
\psi (x)=x^le^{-\gamma x^2+\delta x}\prod\limits_{k=1}^\infty (1+\delta_kx)e^{\delta_kx}
$$
where $l\in\ZZ_+$, $\gamma , \delta , \delta_k\in\RR$, $0<\gamma
^2+\sum\limits_{k=1}^\infty |\delta_k |<\infty$. Then
(see~\cite{MShW},~\cite{BKI}) the kernel
$$
K(x)={\sum\limits_{k=-\infty}^\infty}'\;
\frac{e^ikx}{\psi (ik)}
$$
(${\sum\limits_{k=-\infty}^\infty} '$ denotes that the summation is carried out over all $k$ such that $\psi (ik)\neq 0$) is a $CVD$ -- kernel.
In particular, Bernuolly's kernels
$$
B_r(x)=\frac 1{2\pi}{\sum\limits_{k=-\infty}^\infty}'\;
\frac{e^ikx}{(ik)^r}
$$
and, more generally, kernels of the form
$$
B_{\cal P}(x)=\frac 1{2\pi}{\sum\limits_{k=-\infty}^\infty}'\;
\frac{e^ikx}{{\cal P}(ik)},
$$
where ${\cal P}(x)$ is an algebraic polynomial having real zeros
(integral operator of convolution with such a kernel inverses
differential operators of the form ${\cal P}\left(\frac
d{dx}\right)$) are $CVD$-kernels.

Denote by $F_p \; (1\le p\le \infty)$ the unite ball in the space $L_p$.
 Given a kernel $K$ denote by $K*F_p$ the class of functions $f$ of the form (1)
 where $a\in \RR ,\; \phi\in F_p,\; \phi\perp\mu$.
 Note that $B_r*F_p$ is standard for the theory of quadrature formulas class $W^r_p$
 of real-valued, $2\pi$-periodic functions $f$ having locally absolutely continuous
 derivative $f^{(r-1)}\; (f^{(0)})=f$ and such that $\| f^{(r)}\|_p\le 1$.

Consider the set $Q_n\; (n=1,2,...)$ of all possible quadrature formulas of the form
$$
q(f)=\sum\limits_{k=1}^nc_kf(x_k),\eqno (2)
$$
where $c_k\in \RR , x_1<x_2<...<x_n<x_1+2\pi$.  The problem about optimal for the class $K*F_p$ quadrature formula from $Q_n$ is formulated in the following way.
Find the value
$$
R_n(K*F_p)=\inf\limits_{q\in Q_n}\sup\limits_{f\in K*F_p}\left|\int\limits_0^{2\pi}f(x)dx-q(f) \right|,\eqno(3)
$$
and parameters (knots $x_k$ and coefficients $c_k$) of a quadrature formula $q$ that realizes $\inf$ in the right hand part of (3).

This problem was completely solved for the classes $W^r_p\;
(r=1,2,...; \,1\le p\le \infty )$ in the papers of V. P.
Motornyi~\cite{Motorn}, A. A. Ligun~\cite{Ligun}, A. A.
Zhensykbaev~\cite{Zhensyk},~\cite{Zhensyk1}. It was proved that for
any $n$ the optimal formula has $n$ equidistant knots and equal
coefficients. In the
papers~\cite{Grank},~\cite{BG},~\cite{Chahkiev},~\cite{Chahkiev1}
these results were generalized to the cases of more general function
classes $K*F_p$.

We will consider the following problem. Let $CVD$-kernel $K$; $p\in
[1,\infty ]$; $n=1,2,...$; $h\in (0,\pi /n)$ be given. Denote by
$Q_{n,h}$ the set of all functionals of the form
$$
q(f)=\sum\limits_{k=1}^nc_k\frac
1{2h}\int\limits_{x_k-h}^{x_k+h}f(t)dt\eqno(4)
$$
where $c_k\in \RR$, $x_1+h<x_2-h<x_2+h<...<x_n-h<x_n+h<x_1+2\pi -h$.
Set
$$
R(f,q)=\int\limits_0^{2\pi}f(t)dt-q(f),
$$
$$
R(K*F_p, q)=\sup\limits_{f\in K*F_p}|R(f,q)|,
$$
$$
R_{n,h}(K*F_p)=\inf\limits_{q\in Q_{n,h}}R(K*F_p, q).\eqno(5)
$$
The problem is formulated as follows. Find the value (5) and
parameters $x_k$ and $c_k \; (k=1,...,n)$ of a functional $q$ that
realizes $\inf$ in the right hand part of (5).

From the applications point of view, interval quadrature formulae
are more natural than the usual quadrature formulae based on values
at points, since quite often the result of measuring physical
quantities, due to the structure of the measurement devices, is an
average values of the function, describing the studied quantities,
over some interval. Note that one can obtain the usual quadrature
formula from the corresponding interval quadrature formula as a
limit case, setting $h\to 0$.

In this paper we solve the problem on optimal interval quadrature
formula for the class $K*F_p$ with $p=1$. Thus we will essentially
use Ligun's idea from~\cite{Ligun}.

2. Suppose that $f\in K*F_1$, $q\in Q_{n,h}$, and obtain for
$R(f,q)$ an integral representation. Given $q\in Q_{n,h}$ let
$$
H(q;t)=\frac 1{2h}\sum\limits_{k=1}^nc_k\chi_{(x_k-h,x_k+h)}(t)
$$
where $\chi_A$ is the indicator of a set $A\subset [0, 2\pi]$
continued with the period $2\pi$ on the entire real axis. We will
have
$$
R(f,q)= a\cdot \mu\left( 2\pi
-\sum\limits_{k=1}^nc_k\right)+\int\limits_0^{2\pi}\phi(u)\left[
\int\limits_0^{2\pi}K(t-u)dt-\int\limits_0^{2\pi}K(t-u)H(q;t)dt\right]du
$$
$$
=a\cdot \mu\left( 2\pi
-\sum\limits_{k=1}^nc_k\right)+\int\limits_0^{2\pi}\phi(u)\int\limits_0^{2\pi}K(t-u)\left[
1-H(q;t)\right]dtdu.
$$
Let
$$
M(q,u)=\int\limits_0^{2\pi}K(t-u)\left[ 1-H(q;t)\right]dt.
$$
We obtain
$$
R(f,q)=a\cdot \mu\left( 2\pi
-\sum\limits_{k=1}^nc_k\right)+\int\limits_0^{2\pi}\phi(u)M(q,u)du.\eqno(6)
$$

In the case $\mu (K)=0$ the first term in the right-hand part of (6)
is equal to zero. Solving the problem about optimal interval
quadrature formula for the class $K*F_1$ in the case $\mu (K)=1$ we
can consider functionals $q$ such that $\sum\limits_{k=1}^nc_k=2\pi$
only (otherwise $R(K*F_1,q)=+\infty$). Thus in the case $\mu (K)=1$
we can assume that the first term in the right-hand part of (6) is
equal to zero also.

Taking into account the relations (6), above presented facts, and S.
M. Nikol'skii's duality theorem (see~\cite{Nik},~\cite[Chapt. 2,
Theorem 2.2.1]{Korn}) we obtain
$$
R(K*F_1;q)=\left\{\begin{array}{ll}
                \inf\limits_{\lambda\in\RR}\left\|M(q;\cdot)-\lambda\right\|_\infty ,   &{\rm if}\;\;\;
                 \mu (K)=1, \\ [10pt]
                \left\|M(q;\cdot)\right\|_\infty,  &{\rm if}\;\;\;
                \mu(K)=0.
                 \end{array}\right.\eqno (7)
$$

For $n=1,2,...;\; \lambda\in \RR$, and $0<h<\pi /n$ set
$$
q_{n,h,\lambda}(f)=\lambda \sum\limits_{k=1}^n\frac
1{2h}\int\limits_{2k\pi n^{-1}-h}^{2k\pi n^{-1}+h}f(t)dt.
$$

\begin{theorem}
Let $K$ be a $CVD$-kernel, $n=1,2,...;$, and $0<h<\pi /n$. Then
$$
R_{n,h}(K*F_1)=\left\{\begin{array}{ll}
                R(K*F_1;q_{n,h,2\pi /n})=\inf\limits_{\lambda\in\RR}\left\|M(q_{n,h,2\pi /n};\cdot)-\lambda\right\|_\infty ,   &{\rm if}\;\;\;
                 \mu (K)=1, \\ [10pt]
                R(K*F_1;q_{n,h,\lambda})=\inf\limits_{\lambda\in\RR}\left\|M(q_{n,h,\lambda};\cdot)\right\|_\infty,  &{\rm if}\;\;\;
                \mu(K)=0.
                 \end{array}\right.
$$
\end{theorem}
\medskip

{\bf Proof.} Case $\mu(K)=1$. Denote by $\lambda_q$ ($q\in Q_{n,h}$)
the constant of the best $L_\infty$-approximation of the function
$M(q;\cdot)$. Suppose that for some $q\in Q_{n,h}$ the inequality
$$
\|M(q;\cdot)-\lambda_q\|_\infty <\|M(q_{n,h,2\pi
/n};\cdot)-\lambda_{q_{n,h,2\pi /n}}\|_\infty
$$
holds true. Let
$$
\Delta_\tau (t)=M(q;t-\tau)-M(q_{n,h,2\pi /n};t)-\lambda_q
+\lambda_{q_{n,h,2\pi /n}},\; \tau\in \RR.
$$
Since the function $M(q_{n,h,2\pi /n};\cdot)$ is $2\pi /n$-periodic,
we will have that for any $\tau$
$$
\nu (\Delta_\tau (\cdot))\ge 2n.
$$
It is clear that $K(-\cdot )\in CVD$ with $\mu =1$. Taking into
account this fact and representation (6) it is easy to verify that
for any $\tau$ we will have
$$
\nu (\Delta_\tau '(\cdot))\ge 2n
$$
where
$$
\Delta_\tau ' (t)=-H(q_{n,h,2\pi /n};t)-H(q;t-\tau).
$$

Among the coefficients $c_k$ of the functional $q$ choose such that
$|c_k|\le 2\pi /n$ and denote by $k_0$ its number. Set $\tau_0=2\pi
/n-x_{k_0}$.

\begin{lemma}
$$
\nu (\Delta_{\tau_0}(\cdot))\le 2n-2.
$$
\end{lemma}

If Lemma 1 will be proved, we obtain the contradiction with the
above presented statement about the number of the sign changes of
$\Delta_\tau '$. After that the part of Theorem 1 related to the
case $\mu (K)=1$ will be proved.

\medskip

{\bf Proof of the Lemma 1.} Suppose that $\nu (\Delta_{\tau_{0}}
')\ge 2n$. It means that there exist points
$y_1<y_2<...<y_{2n}<y_1+2\pi$ such that difference
$\Delta_{\tau_{0}}'$ has at these points nonzero values of
alternating sign. Let (for definiteness) $\Delta_{\tau_{0}}'$ has
negative values at the points $y_2,y_4,...,y_{2n}$. It is easily
seen that such points  must belong to the different intervals of
positivity of the function $H(q;t-\tau_0)$. But the number of such
intervals inside an interval of the length $2\pi $ is less than or
equal to $n$. Really, the difference $\Delta_{\tau_{0}}'$ can not
change sign inside the common part of an interval of positivity of
the function $H(q_{n,h,2\pi /n};t)$ and an interval of positivity of
the function $H(q;t-\tau_0)$ (this difference is constant inside
such an part). Thus at least one of the point of negativity (say
$y_{2l}$) must belongs to the interval $(2\pi /n-h, 2\pi /n+h)$. But
in view of the choice of $\tau_0$, the difference
$\Delta_{\tau_{0}}'$ is nonnegative on this interval. This is a
contradiction.

Lemma is proved.

\medskip

The case $\mu (K)=0$. We need the following analog of the Lemma 1
from~\cite{BG}.

\begin{lemma}
Let $\lambda ^*$ be such that
$$
\inf\limits_{\lambda\in \RR}\| M(q_{n,h,\lambda};\cdot)\|_\infty =\|
M(q_{n,h,\lambda ^*};\cdot)\|_\infty.
$$
Then
$$
\max\limits_t M(q_{n,h,\lambda ^*};t)=-\min\limits_t
M(q_{n,h,\lambda ^*};t).
$$
\end{lemma}

\medskip

{\bf Proof of Lemma 2.} We have
$$
M(q_{n,h,\lambda
^*};t)=\int\limits_0^{2\pi}K(u)du-\lambda\int\limits_0^{2\pi}K(u-t)H(q_{n,h,1};u)du.
$$
It is well known that $CVD$-kernel $K$ with $\mu (K)=0$ does not
change  sign (let it is nonnegative for definiteness). Then it is
easily seen that function
$$
\psi (t)=\int\limits_0^{2\pi}K(u-t)H(q_{n,h,1};u)du
$$
is strictly positive. If for a given $\lambda_0>0$
$$
\max\limits_t\left[ \int\limits_0^{2\pi}K(u)du-\lambda_0\psi
(t)\right]>-\min\limits_t\left[
\int\limits_0^{2\pi}K(u)du-\lambda_0\psi (t)\right],\eqno(8)
$$
then for $\lambda <\lambda_0$ and close enough to $\lambda_0$ we
will have
$$
\left\| \int\limits_0^{2\pi}K(u)du-\lambda\psi (\cdot
)\right\|_\infty < \left\| \int\limits_0^{2\pi}K(u)du-\lambda_0\psi
(\cdot) \right\|_\infty.\eqno (9)
$$
If instead of (8) we have
$$
\max\limits_t\left[ \int\limits_0^{2\pi}K(u)du-\lambda_0\psi
(t)\right]<-\min\limits_t\left[
\int\limits_0^{2\pi}K(u)du-\lambda_0\psi (t)\right],
$$
then the inequality (9) will hold true for all $\lambda >\lambda_0$
and close enough to $\lambda_0$. Since the existence of $\lambda^*$
and its positiveness are obvious, the Lemma is proved.

\medskip

Suppose that for some $q\in Q_{n,h}$
$$
R(K*F_1; q)<\inf\limits_{\lambda\in \RR}R(K*F_1;
q_{n,h,\lambda})\eqno (10)
$$
or that is equivalent (in view of (7))
$$
\| M(q;\cdot)\|_\infty<\|
M(q_{n,h,\lambda^*};\cdot)\|_\infty.\eqno(11)
$$
Taking into account the Lemma 2 and the fact that
$M(q_{n,h,\lambda^*};\cdot)$ is $2\pi /n$-periodic, we conclude that
$$
\nu (\Delta_\tau (\cdot))\ge 2n.
$$
where
$$
\Delta _\tau (t)=M(q;t-\tau)-M(q_{n,h,\lambda^*};t)
$$
But
$$
\Delta _\tau
(t)=\int\limits_0^{2\pi}K(u-t)[H(q_{n,h,\lambda^*};u)-H(q;u-\tau)]du.
$$
Since $K\in CVD$ for any $\tau$ then
$$
\nu(\Delta_\tau '):=\nu (H(q_{n,h,\lambda^*};u)-H(q;u-\tau))\ge 2n.
$$
However analogously to the case $\mu (K)=1$ it is possible to choose
$\tau_0$ such that $\nu(\Delta_{\tau_0} ')\le 2n-2$. We omit the
details.

Therefore the relation (11), and consequently the relation (10), is
impossible. Theorem 1 is proved.

\bigskip

Let us show that Theorem 1 implies the optimality of the quadrature
formula with equidistant knots and equal coefficients on the class
$K*F_1$ among the all quadrature formulas from $Q_n$. Let a
$CVD$-kernel $K$ be such that the class $K*F_1$ is a relatively
compact subset of the space $C$ (in the case $\mu (K)=0$) or can be
obtained by shifts to the constants from relatively compact subset
of $C$ (in the case $\mu (K)=1$). For any quadrature formula $q$ of
the form (2) let $\{q_h,\;h>0\}$ be  we the family of functionals of
the form (4) having the same $c_k$ and $x_k$. It is easily seen that
$q_h(f)\to q(f)$ uniformly on the set $f\in K*F_1$ as $h\to 0$
(remind that in the case $\mu (K)=1$ we consider formulas $q$ and
functionals $q_h$ such that $\sum\limits_{k=1}^nc_k=2\pi$ only). Set
$$
q_{n,0,\lambda}(f)=\lambda\sum\limits_{k=1}^nf\left( \frac {2k\pi
}{n}\right).
$$
Then (we restrict ourself by the case $\mu (K)=1$) for any $q\in
Q_n$ we will have
$$
R(K*F_1;q)=\lim\limits_{h\to 0}R(K*F_1;q_h)\ge\lim\limits_{h\to
0}R(K*F_1;q_{n,h,2\pi /n})=R(K*F_1;q_{n,0,2\pi /n})
$$
that is the optimality of the rectangle formula on the class
$K*F_1$.

\medskip

Finely we note that the statement of the Theorem 1 can be easily
generalized to the classes of convolutions of the functions from
$F_1$ with $O(M, \Delta )$-kernels (see~\cite{BG}).

\end{document}